\documentclass{article}
\usepackage{arxiv}
\usepackage{lineno,hyperref}
\modulolinenumbers[5]
\usepackage{multirow}
\usepackage{algorithm}
\usepackage{algpseudocode}

\usepackage{amssymb}
\usepackage{latexsym}
\usepackage{bm}

\usepackage{url}
\usepackage{xcolor}
\usepackage{graphicx}
\usepackage{rotating}
\usepackage{adjustbox}

\usepackage{epstopdf}
\usepackage{subcaption}
\usepackage{amsmath}
\usepackage{amsthm}

\usepackage{tikz}

\DeclareRobustCommand\full  {\tikz[baseline=-0.6ex]\draw[thick] (0,0)--(0.5,0);}
\DeclareRobustCommand\dotted{\tikz[baseline=-0.6ex]\draw[thick,dotted] (0,0)--(0.54,0);}
\DeclareRobustCommand\dashed{\tikz[baseline=-0.6ex]\draw[thick,dashed] (0,0)--(0.54,0);}

\DeclareRobustCommand\dashdot {\tikz[baseline=-0.6ex]\draw[thick,dash dot] (0,0)--(0.5,0);}

\algnewcommand\algorithmicinput{\textbf{Input:}}
\algnewcommand\INPUT{\item[\algorithmicinput]}
\algnewcommand\algorithmicoutput{\textbf{Output:}}
\algnewcommand\OUTPUT{\item[\algorithmicoutput]}
\algnewcommand\PARAM{\item[\textbf{Parameters:}]}

\newtheorem{theorem}{Theorem}

\title{Monomial barrier functions for the box-constrained convex optimization problems}

\author{ Hatem A.~Fayed \\
	University of Science and Technology, Mathematics Program, Zewail City of Science and Technology\\
October Gardens, 6th of
October, Giza 12578, Egypt\\
	\texttt{hfayed@zewailcity.edu.eg} \\
}

\renewcommand{\headeright}{}
\renewcommand{\undertitle}{}

\hypersetup{
pdftitle={Monomial barrier for box-constraints},
pdfsubject={Monomial barrier for box-constraints},
pdfauthor={Hatem A.~Fayed},
pdfkeywords={Convex optimization, Box-constraints, Newton's method, Barrier Methods},
}

\begin{document}
\maketitle

\begin{abstract}
In this article, a novel barrier function is introduced to convert the box-constrained convex optimization problem to an unconstrained problem. For each double-sided bounded variable, a single monomial function is added as a barrier function to the objective function. This function has the properties of being positive, approaching zero for the interior/boundary points and becomes very large for the exterior points as the penalty parameter approaches zero. The unconstrained problem can be solved efficiently using Newton's method with a backtracking line search. Experiments were conducted using the proposed method, the interior-point for the logarithmic barrier (IP), the trust-region reflective (TR) and the limited-memory Broyden, Fletcher, Goldfarb, and Shanno for bound constrained problems (LBFGSB) methods on the convex quadratic problems of the CUTEst collection. Although the proposed method was implemented in MATLAB, the results showed that it outperformed IP and TR for all problems. The results also showed that despite LBFGSB was the fastest method for many problems, it failed to converge to the optimal solution for some problems and took a very long time to terminate. On the other hand, the proposed method was the fastest method for such problems. Moreover, the proposed method has other advantages, such as: it is very simple and can be easily implemented and its performance is expected to be improved if it is implemented using a low-level language, such as C++ or FORTRAN on a GPU.

 \end{abstract}

\section{Introduction}

The general simple box-constrained convex optimization problem considered in this article is,

\begin{equation}
\label{eq1}
\min _{\boldsymbol{x} \in\Omega} \,f(\boldsymbol{x}),\quad\   s.t. \,\Omega =\left\{\boldsymbol{x} \in \mathbb{R}^{m}: \boldsymbol{l} \leq \boldsymbol{x} \leq \boldsymbol{u}\right\}
\end{equation}

where $\boldsymbol{x},\boldsymbol{l},\boldsymbol{u}\in \mathbb{R}^m$ and $f(\boldsymbol{x}) $ is twice continuously differentiable convex function.

The Karush–Kuhn–Tucker (KKT) conditions at the optimal solution, $\boldsymbol{x}^*$, are given by \cite{Andrei17},
\begin{equation}
\label{eq1opt}
\begin{aligned}
x^*_i=l_i \Rightarrow &\nabla f_i (\boldsymbol{x}^*)\ge 0   \\
l_i<x^*_i<u_i \Rightarrow &\nabla f_i (\boldsymbol{x}^*)= 0   \\
x^*_i=u_i \Rightarrow &\nabla f_i (\boldsymbol{x}^*) \le 0   \\
\end{aligned}
\end{equation}
for $i=1\cdots m$ where $\nabla f$ is the gradient of $f$.

This type of problems arises in many applications such as: the constrained autoencoders feature extraction for deep learning \cite{Ayinde18}, the water flow through a porous dam \cite{Baiocchi73}, the optimal design problem \cite{Birgin99}, the torsion applied to a bar \cite{Cea73}, the journal bearing \cite{Cryer71}, the simulation of time-dependent contact and friction problems in rigid body mechanics \cite{Lotstedt84}, the obstacle problem \cite{Rodrigues20} and the linear least-squares problems with upper and lower bounds \cite{Stark95}.

There are many algorithms in the literature were proposed to solve this problem, some are restricted to convex functions (in some cases quadratic) and some are general. For the convex quadratic function, the problem can be solved in polynomial time by interior-point or trust-region algorithms \cite{Goldfarb90,Heinkenschloss99}, however, for an indefinite quadratic function, the problem is NP-hard \cite{Murty87}. In the sequel, the methods used for the convex objective functions are overviewed. 

Early methods include the active set method \cite{Fletcher74,Gill81,OLeary80}, in which a sequence of sub-problems are solved iteratively. Initially, it identifies some variables as active (fixed at the bounds) defining a face likely containing the stationary point and then the objective function is minimized using an unconstrained algorithm with respect to the remaining (inactive) variables. In the successive iterations, one constraint is dropped/added to the active set. For large-scale problems, this would be impractical as many iterations would be required to achieve the correct active set. Later, several variants were proposed for the convex quadratic function that are based on the projected Newton method which can add and drop many constraints at each iteration \cite{Bertsekas82,Dembo83,Yang91}. Several research and convergence study for the convex quadratic function were explored in \cite{Dostal97,Dostal03,More89,More91}. 

Bertsekas \cite{Bertsekas82a} applied the gradient projection to the active set algorithm for the general convex function. Further extensions and study of the convergence can be found in \cite{Bertsekas82a,Calamai87}. Byrd et al. \cite{Byrd95} developed the LBFGSB algorithm which approximates the Hessian of the objective function and performs the gradient projection method to determine the active set. Then, it performs line searches along the directions obtained by a limited memory BFGS method \cite{Byrd94} to explore the subspace of the inactive variables. An implementation of their algorithm is available in~\cite{Zhu97}. Ni and Yuan \cite{Ni97} used the limited memory quasi-Newton methods to update the inactive variables and a projected gradient method to update the active variables. The computational cost and memory requirements per each iteration are small for these methods. This leads to an efficient implementation for large problems. 

Lin and Moré \cite{Lin99} used a trust-region Newton's method which relies on the geometry of the feasible set that generates a Cauchy step, a preconditioned conjugate gradient method with an incomplete Cholesky decomposition to generate a direction, and a projected search to compute the step. Coleman and Li \cite{Coleman94,Coleman96a} suggested a different approach that is related to the trust-region method, namely, the affine-scaling interior-point method. In this method, the first order optimality conditions of (\ref{eq1}) is rewritten as a nonlinear system of equations using a certain scaling matrix. The resulting system is solved using Newton's method. Some developments of this method were presented in \cite{Heinkenschloss99,Ulbrich99}.

Recently, Hager and Zhang \cite{Hager06} developed an active set algorithm (ASA) for box constrained optimization which consists of a nonmonotone gradient projection step, an unconstrained optimization step, and a set of rules for branching between the two steps. In \cite{Xu07}, ASTRAL algorithm, an active-set algorithm that uses both the active-set identification and LBFGS \cite{Liu89} updating for the Hessian approximation, was developed. It uses a gradient projection step to determine an active face and then solves trust-region subproblems using primal-dual interior-point techniques. In \cite{Cristofari17}, a two-satge active set algorithm was proposed. In this algorithm, an active set is first estimated by applying a non-monotone line search algorithm, then a truncated-Newton strategy is used in the subspace of the non-active variables. 

\section{The monomial barrier function}\label{sec4}
\begin{figure}
    \centering
     \includegraphics[width=0.6\textwidth]{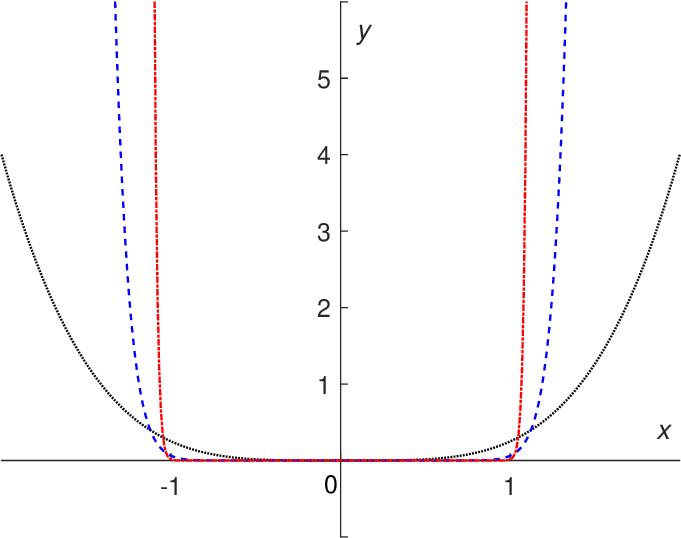}
    \caption{Graphs of $x^\mu/\mu$ for different values of $\mu$. $\mu=4$ (\dotted) , $\mu=16$ (\textcolor{blue}\dashed) and $\mu=64$ (\textcolor{red}\dashdot).}
    \label{monomial_Fig}
\end{figure}

The main idea of the proposed monomial barrier function stems from the shape of the function $g\left(x\right)= x^\mu/\mu$ where $\mu$ is an even positive integer number (see Fig. \ref{monomial_Fig}). For a sufficiently large $\mu$, this function tends to be zero for $\lvert x \rvert \le 1$ and very large for $\lvert x \rvert >1$. So, it can be used as a barrier function for the box constraint: $-1\le x \le 1$. 

In the following example, a simple 1-dimensional problem is introduced to illustrate the different cases where the minimizer is an interior, exterior and boundary point respectively.

{\bf Example 1}
\label{ex1}
Consider the following optimization problem,
\begin{equation}
\label{eq4}
\begin{aligned}
\min _{x} f(x) &= \frac{1}{2} \left(x-a\right)^2 \quad
s.t. &-1 \leq x \leq 1
\end{aligned}
\end{equation}

The monomial barrier function can be defined as,

\begin{equation}
\label{eq5}
P(x;\mu) = \frac{1}{2} \left(x-a\right)^2+x^\mu/\mu
\end{equation}

As $\mu\to\infty$ the minimizer of $P(x;\mu)$ approaches the solution of (\ref{eq4}).
Note that $P(x;\mu)$ is always nonnegative. For a certain $\mu$, the exact minimizer of $P(x;\mu)$ can be obtained by setting $\frac{\partial P(x;\mu)}{\partial x}=0$; i.e. 
\begin{equation}
\label{eq7}
x-a+x^{\mu-1}=0
\end{equation}

which can be viewed as the point of intersection of the following curves,

\begin{equation}
\label{eq8}
y=x^{\mu-1},y=a-x
\end{equation}

\begin{figure}[h!]
    \centering 
\begin{subfigure}{0.45\textwidth}
  \includegraphics[width=\linewidth]{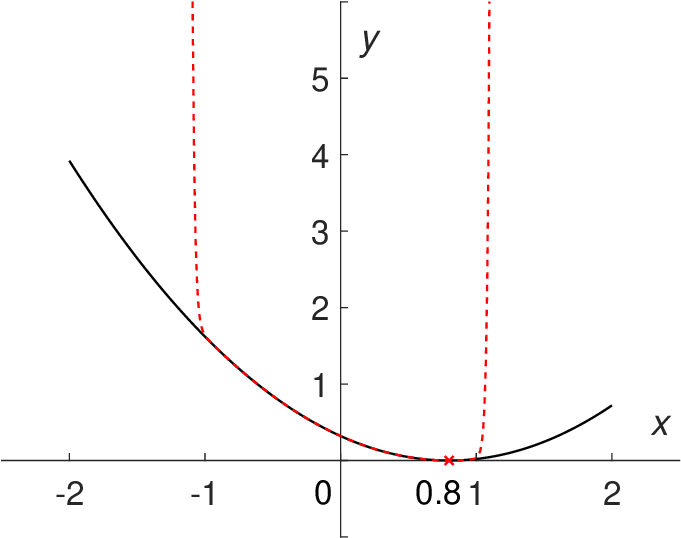}
  \caption{$a=0.8$.}
  \label{fig:1}
\end{subfigure}
\hfill
\begin{subfigure}{0.45\textwidth}
  \includegraphics[width=\linewidth]{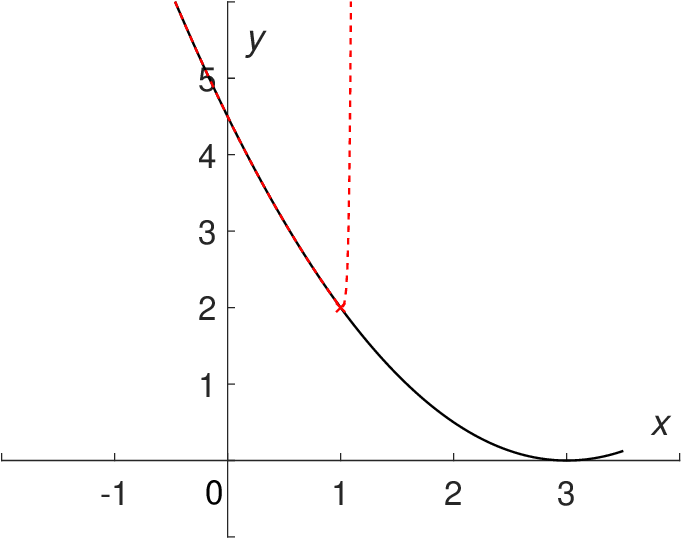}
  \caption{$a=3$.}
  \label{fig:2}
\end{subfigure}
\hfill
\begin{subfigure}{0.45\textwidth}
  \includegraphics[width=\linewidth]{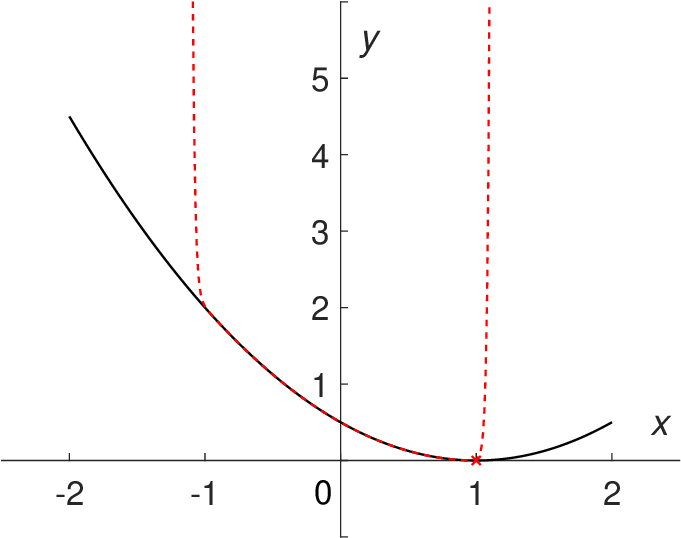}
  \caption{$a=1$.}
  \label{fig:3}
\end{subfigure}

\caption{Graphs of $f(x)$ (\full) and $P(x;64)$ (\textcolor{red}\dashed) for Example 1 for different values of $a$.}
\label{example1_Fig}
\end{figure}

Three different cases  can be obtained according to the value of $a$ (see Fig. \ref{example1_Fig}).

\textit{Case I:} $-1<a<1$

\begin{figure}[ht]
\centering
\includegraphics[width=0.6\textwidth]{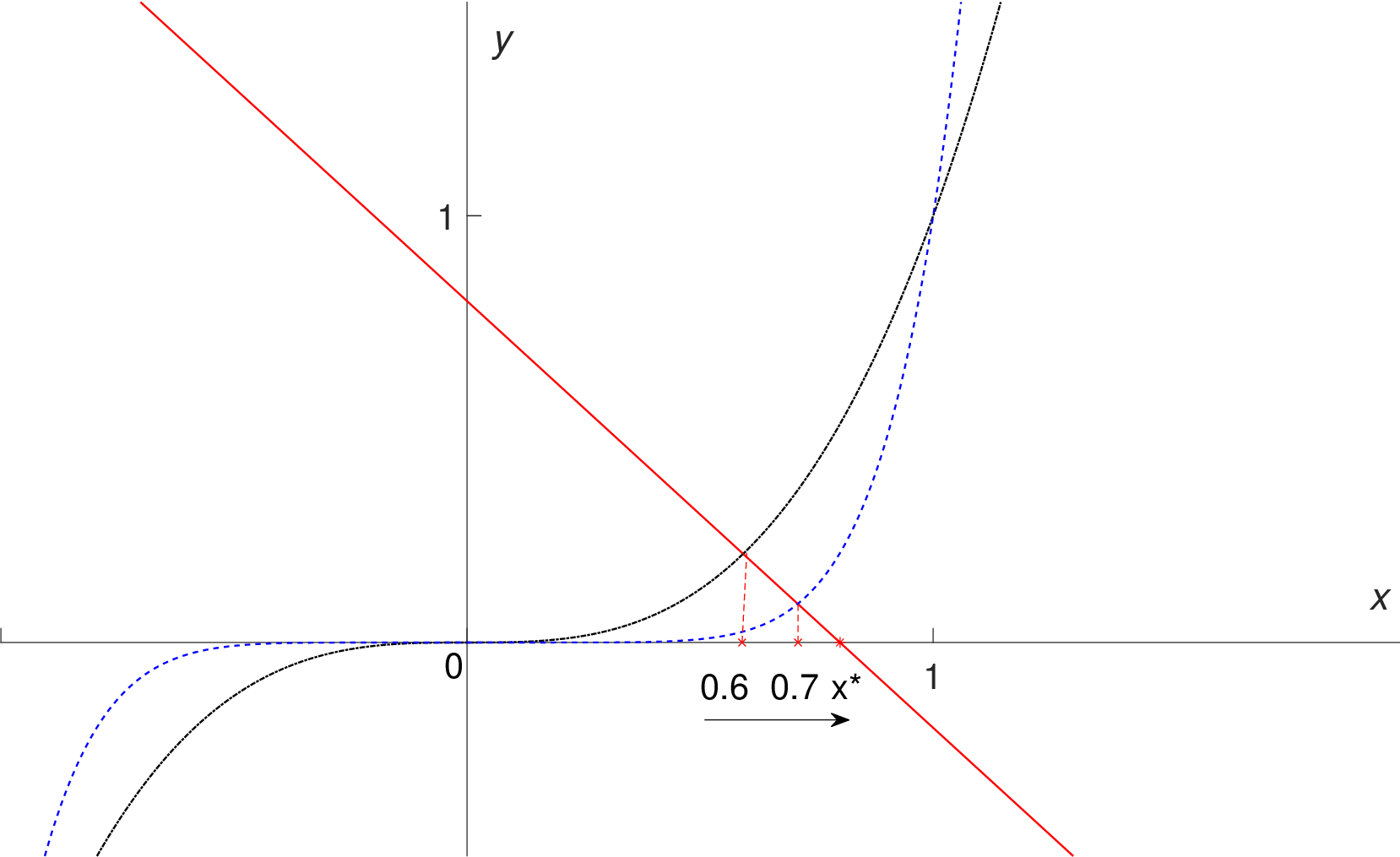}
\caption{Graphs of $y=0.8-x$(\textcolor{red}\full), $y=x^3$ (\dashdot) and $y=x^7$ (\textcolor{blue} \dashed).
}\label{inter_Fig}
\end{figure}

Clearly, the optimal solution is $x^*=a$ and the optimal function value is $f^*=0$. Consider $a=0.8$, for $\mu=4$, the solution of (\ref{eq7}) is $x\approx0.6$ while for $\mu=8$, the solution is $x\approx0.7$ which is closer to the optimal solution, $x^*=0.8$ (see Fig. \ref{inter_Fig}). For $\mu=1024$, the solution of (\ref{eq7}), approximated to $8$ decimal places, is: $x=0.80000000$. As expected, the higher the value of $\mu$, the better the approximate solution. Graphs of $f(x)$ and $P(x;64)$ are shown in Fig. \ref{example1_Fig}\subref{fig:1}.

\textit{Case II:} $\lvert a \rvert>1$

Consider $a=3$, in this case, the optimal solution is $x^*=1$ and the optimal function value is $f^*=2$. For $\mu=4$, the solution of (\ref{eq7}) is $x\approx1.2$ while for $\mu=8$, the solution is $x\approx1.1$. Although both are infeasible, their projections onto the feasible region lead to the optimal solution, $x^*=1$. This can be incorporated to achieve a fast convergence as will be described later in algorithm \ref{algo1}. Graphs of $f(x)$ and $P(x;64)$ are shown in Fig. \ref{example1_Fig}\subref{fig:2}.

\textit{Case III:} $a=\pm 1$

Consider $a=1$, in this case, the optimal solution is $x^*=1$ and the optimal function value is $f^*=0$. For $\mu=1024$, the solution of (\ref{eq7}), approximated to $8$ decimal places, is: $x=0.99486088$. Thus, for the boundary points, larger values of $\mu$ have to be considered to achieve a good approximation. Graphs of $f(x)$ and $P(x;64)$ are shown in Fig. \ref{example1_Fig}\subref{fig:3}.

{\bf Example 2}
Consider the following 2-dimensional optimization problem,
\begin{equation}
\label{eq9}
\begin{aligned}
\min _{\boldsymbol{x}} f(\boldsymbol{x}) &= \frac{1}{2} \left[\left(x_1+2\right)^2+x_2^2\right]  \\
s.t. &-1 \leq x_1 \leq 1,-1 \leq x_2 \leq 1
\end{aligned}
\end{equation}

It is easy to note that the optimal function value occurs at $x_1=-1, x_2=0$. The monomial barrier function for this problem can be defined as,

\begin{equation}
\label{eq10}
P(\boldsymbol{x};\mu) = \frac{1}{2} \left[\left(x_1+2\right)^2+x_2^2\right]+\frac{1}{2 \mu}\left[x_1^\mu+x_2^\mu\right]
\end{equation}

The proposed monomial barrier function can be compared with the logarithmic barrier function which can be defined for this problem as,
\begin{equation}
\label{eq11}
\begin{aligned}
P_l(\boldsymbol{x};\nu) &=\frac{1}{2} \left[\left(x_1+2\right)^2+x_2^2\right]  \\
&-\nu\left[log(1+x_1)+log(1-x_1)+log(1+x_2)+log(1-x_2)\right]
\end{aligned}
\end{equation}
A sequence of values $\{\nu_k\}$ with $\nu_k\to0$ is incorporated, and the approximate minimizer, $\boldsymbol{x}^k$, of $P_l(\boldsymbol{x};\nu_k)$ is obtained for each $k$.

\begin{figure}[ht]
    \centering 
\begin{subfigure}{0.4\textwidth}
  \includegraphics[width=\linewidth]{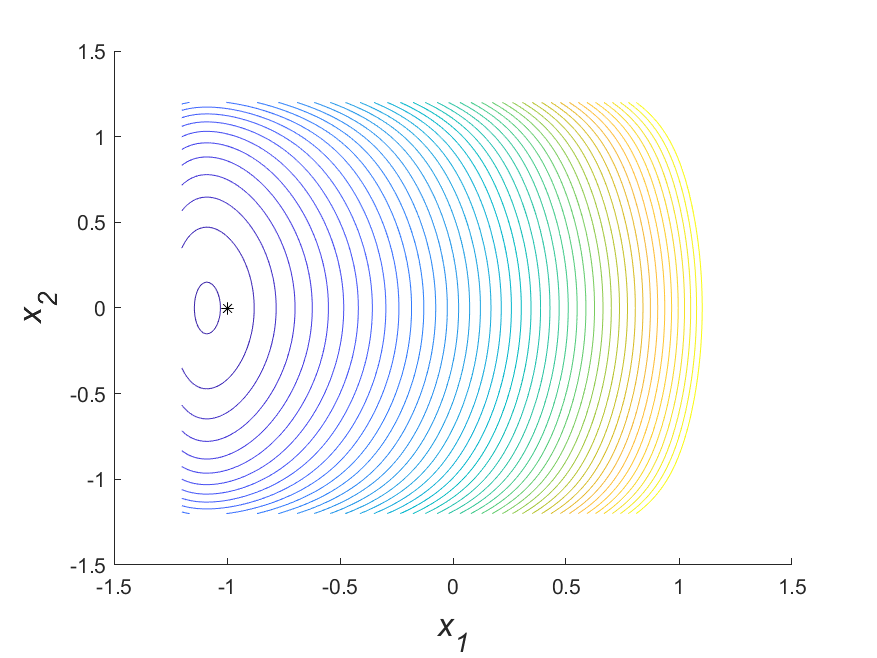}
  \caption{$P(\boldsymbol{x};2^3)$.}
  \label{fig:1}
\end{subfigure}
\hfill
\begin{subfigure}{0.4\textwidth}
  \includegraphics[width=\linewidth]{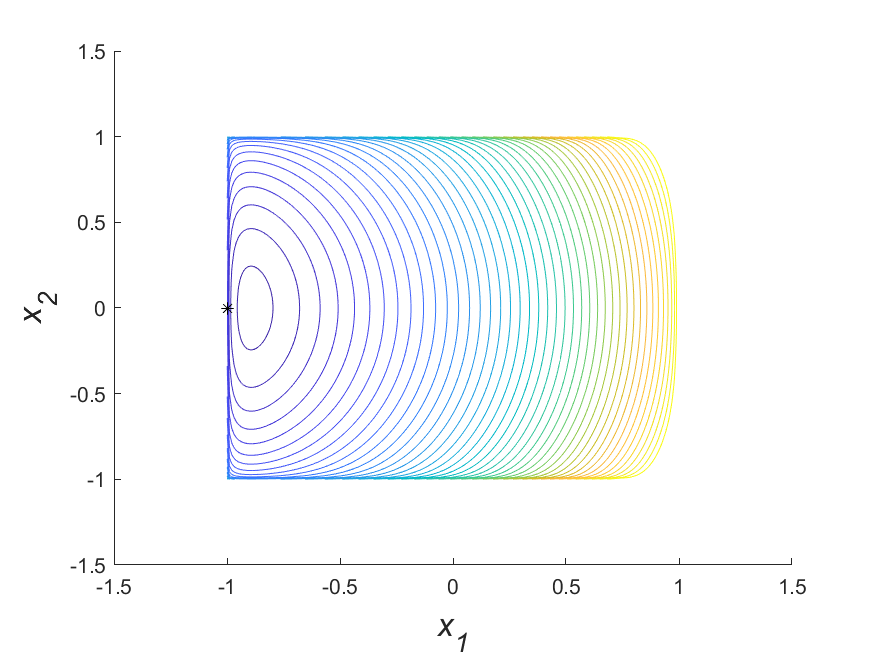}
  \caption{$P_l(\boldsymbol{x};1/2^3)$.}
  \label{fig:2}
\end{subfigure}

\begin{subfigure}{0.4\textwidth}
  \includegraphics[width=\linewidth]{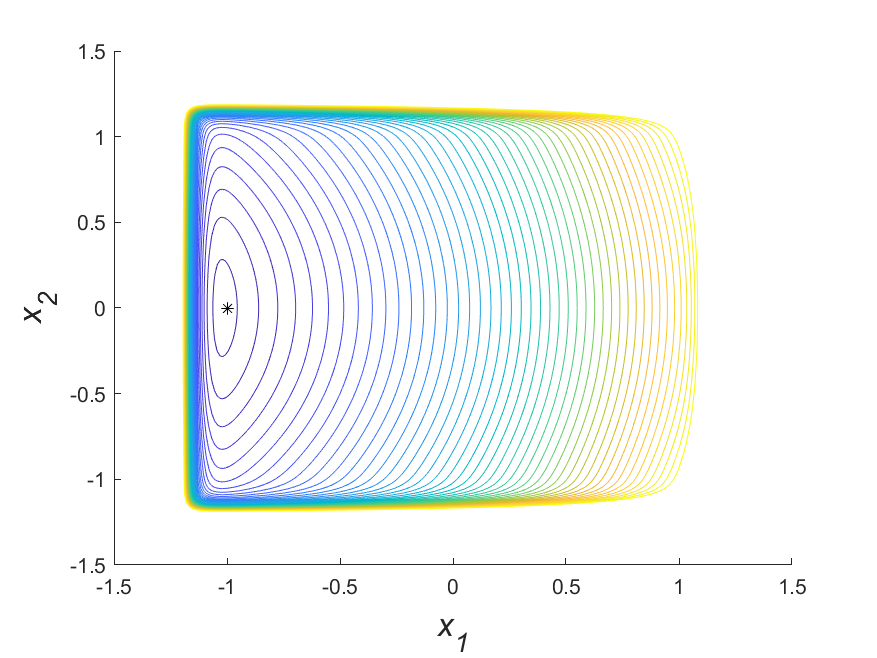}
  \caption{$P(\boldsymbol{x};2^{5})$.}
  \label{fig:3}
\end{subfigure}
\hfill
\begin{subfigure}{0.4\textwidth}
  \includegraphics[width=\linewidth]{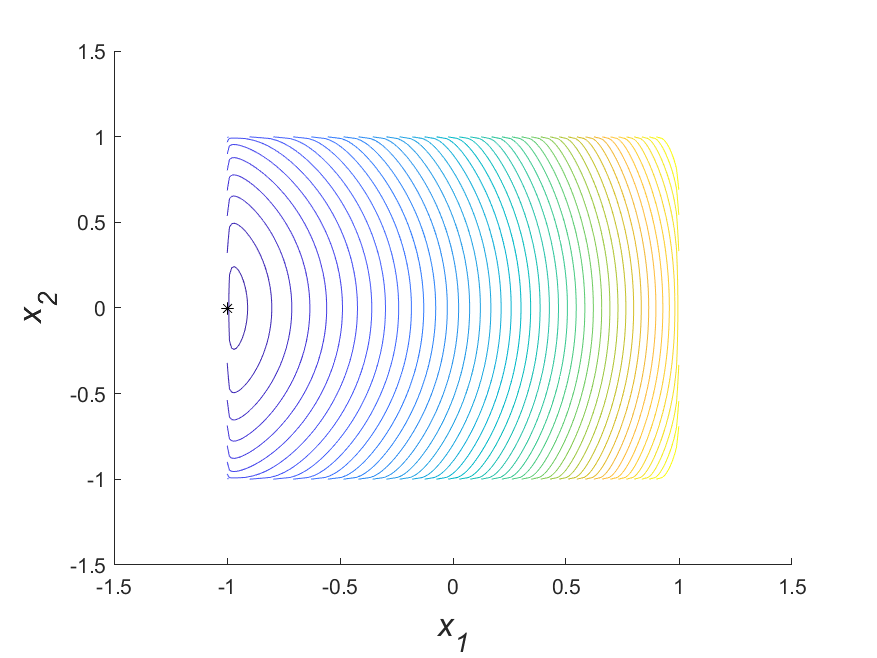}
  \caption{$P_l(\boldsymbol{x};1/2^{5})$.}
  \label{fig:4}
\end{subfigure}

\begin{subfigure}{0.4\textwidth}
  \includegraphics[width=\linewidth]{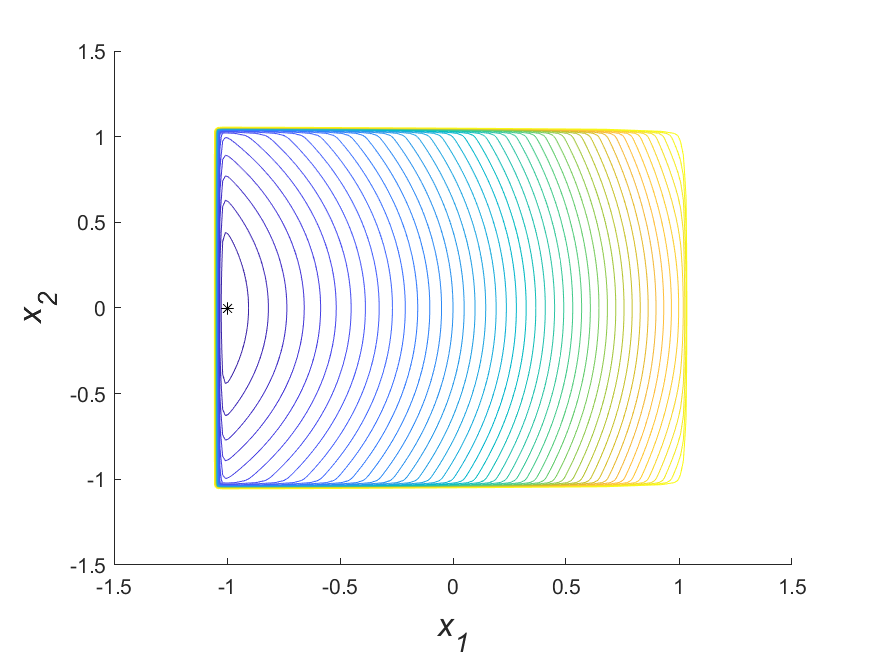}
  \caption{$P(\boldsymbol{x};2^{7})$.}
  \label{fig:5}
\end{subfigure}
\hfill
\begin{subfigure}{0.4\textwidth}
  \includegraphics[width=\linewidth]{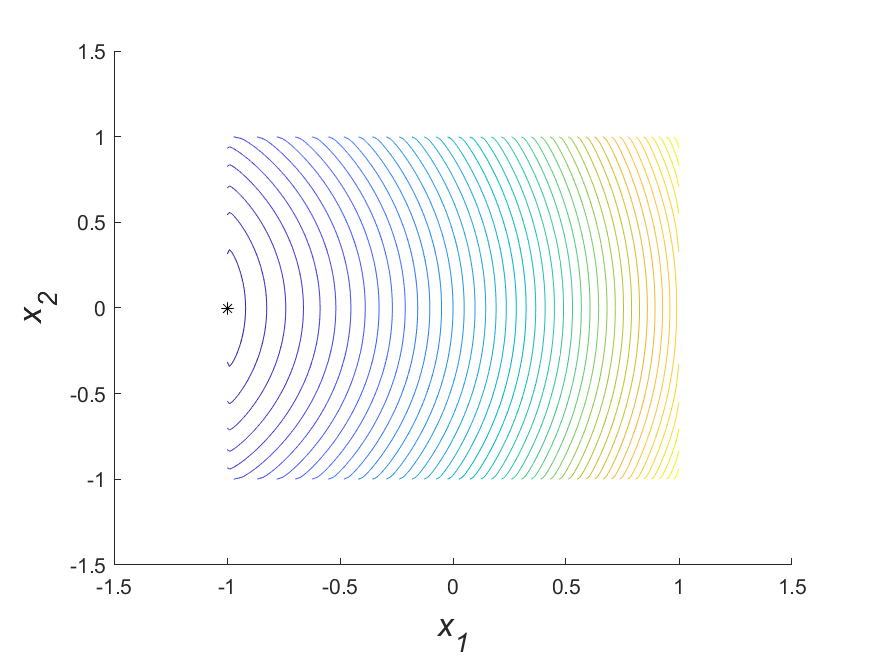}
  \caption{$P_l(\boldsymbol{x};1/2^{7})$.}
  \label{fig:6}
\end{subfigure}

\caption{Contours of $P(\boldsymbol{x};\mu)$ and $P_l(\boldsymbol{x};1/\mu)$ for Example 2 for different values of $\mu$.}
\label{example2_Fig}
\end{figure}

Contours of both $P(\boldsymbol{x};\mu)$ and $P_l(\boldsymbol{x};1/\mu)$ are shown in Fig. \ref{example2_Fig}. It can be noticed that, for $\mu=2^3$, the contours of both of them are elliptical that are not too elongated and hence most unconstrained minimization algorithms can be applied successfully to find the minimizers. The minimizers of $P(\boldsymbol{x};\mu)$ and $P_l(\boldsymbol{x};1/\mu)$ are (-1.00000003,0) and (-0.89340298,0) respectively. While for $\mu=2^{10}$, the minimizers of are (-1.00000000,0) and (-0.99902500,0) respectively. Moreover, the projection of the optimal solution of $P(\boldsymbol{x};\mu)$ onto the feasible region leads directly to the optimal solution of (\ref{eq9}).

It can be noted that, for an active constraint, the contours of $P_l(\boldsymbol{x};\nu)$ are too elongated and its Hessian becomes increasingly ill-conditioned near the minimizer as $\nu \to 0$. This ill-conditioning usually causes problems when solving the linear equations to calculate the Newton step, so, the modified Cholesky decomposition is often used to avoid this problem. On the other hand, the contours of $P(\boldsymbol{x};\mu)$ are not elongated that much near the optimal point, as $\boldsymbol{x}$ is allowed to be feasible or infeasible. This suggests that Cholesky decomposition can be used rather than the modified Cholesky decomposition if $f(\boldsymbol{x})$ is strictly feasible. Moreover, the required system to be solved at each iteration for the monomial barrier method is much smaller than that of the logarithmic barrier method. Hence, the implementation of the monomial barrier method is expected to be faster than the interior-point for the logarithmic barrier method. 

\subsection{General framework}\label{GF}
There are two approaches to handle the general boundary box constraints in (\ref{eq1}). The first one is to scale the problem using the transformations,
\begin{equation}
\label{eqtx}
z_i=\frac{x_{i}-r_{i}}{q_{i}} , i=1\cdots m
\end{equation}
where
\begin{equation}
\label{eqtx3}
 r_{i}=\frac{u_{i}+l_{i}}{2},
 q_{i}=\frac{u_{i}-l_{i}}{2} 
\end{equation}
Then, the problem is transformed into,
\begin{equation}
\label{eq100}
\min _{\boldsymbol{z} \in\Omega} \,f(\boldsymbol{z}),\quad\   s.t. \,\Omega =\left\{\boldsymbol{z} \in \mathbb{R}^{m}: -\boldsymbol{1} \leq \boldsymbol{z} \leq \boldsymbol{1}\right\}
\end{equation}
where $\boldsymbol{1}$ is a vector of ones. Hence the following barrier function can be used,

\begin{equation}
\label{eq101}
P(\boldsymbol{z};\mu) =f(\boldsymbol{z})+\frac{1}{m \mu} \sum_{i=1}^{m}z_i^\mu
\end{equation}
Imposing $1/m$ in the barrier term is used to let the number of outer iterations independent on the problem size as will be shown in Theorem \ref{thm2}.
 
Alternatively, the following barrier function for (\ref{eq1}) can be directly defined as,
\begin{equation}
\label{eq12}
P(\boldsymbol{x};\mu) =f(\boldsymbol{x})+\frac{1}{m \mu} \sum_{i=1}^{m}\left(\frac{x_{i}-r_{i}}{q_{i}}\right)^{\mu}
\end{equation}

The main steps of the monomial barrier algorithm are shown in Algorithm \ref{algo1}. A sequence of values $\{\mu_k\}$ with $\mu_k\to\mu_{max}$ is incorporated, these are called the outer iterations. In each outer iteration, an approximate minimizer $\boldsymbol{x}^k$ of $P(\boldsymbol{x};\mu_k)$ is obtained using Newton's method with a backtracking line search, these are called the inner iterations. This minimizer is then projected onto the feasible region using (\ref{eq_proj}),
\begin{equation}
\label{eq_proj}
proj(x_i,l_i,u_i)= \begin{cases}l_{i} & \text { if } x_{i}<l_{i} \\ x_{i} & \text { if } x_{i} \in\left[l_{i}, u_{i}\right] \\ u_{i} & \text { if } x_{i}>u_{i}\end{cases}
\end{equation}
 
Therefore, a new value $\mu_{k+1}=\tau \mu_k$ is set for the next outer iteration where $\tau>1$ is a positive integer number and the starting point of the next iteration is set to $\boldsymbol{x}^k$. This procedure continues until a required tolerance, $\epsilon_g$, of the projected gradient given by (\ref{eq_gr}) is achieved or $\mu$ exceeds an upper bound, $\mu_{max}$.

\begin{equation}
\label{eq_gr}
\Vert proj\left(\boldsymbol{x}_k-\boldsymbol{g}_k,\boldsymbol{l},\boldsymbol{u}\right) - \boldsymbol{x}_k \Vert_\infty < \epsilon_g
\end{equation}

\begin{algorithm}[ht]
\caption{Monomial Barrier}\label{algo1}
\begin{algorithmic}[1]
\Require $\mu_0$: an even positive integer, $\mu_{max}$: an upper bound for $\mu$, $\tau>1$: a positive integer number, $\boldsymbol{x_0}$: starting point.  
\Ensure $\boldsymbol{x}^*$ 
\For{$k =0,1,\cdots$}
	   \State Find an approximate minimizer $\boldsymbol{x}^k$ for $P(\boldsymbol{x};\mu_k)$ starting at $\boldsymbol{x}_0$ and terminate when the stopping criteria for the inner iterations are met	   
	    \State $\boldsymbol{x}^{k}\leftarrow$ projection of $\boldsymbol{x}^{k}$ onto the bound constraints $\boldsymbol{l} \le \boldsymbol{x}\le \boldsymbol{u}$ 
	    \State $\mu_{k+1} \leftarrow \tau \mu_k$
		\If{ the stopping criteria for the outer iterations are met}
		        \State Terminate with an approximate solution $\boldsymbol{x}^* \leftarrow \boldsymbol{x}^{k}$
		\EndIf
	    \State $\boldsymbol{x}_0 \leftarrow \boldsymbol{x}^{k}$
\EndFor
\end{algorithmic}
\end{algorithm}

For the inner iterations, Newton's method are applied iteratively to find the minimizer for $P(\boldsymbol{x};\mu_k)$. To achieve this, the gradient and Hessian matrix of $P(\boldsymbol{x};\mu)$ are calculated as follows,
\begin{equation}
\label{eq13}
\begin{aligned}
\nabla P(\boldsymbol{x};\mu) &= \nabla f(\boldsymbol{x})+\boldsymbol{e}(\boldsymbol{x};\mu), \\
\nabla^2 P(\boldsymbol{x};\mu) &= \nabla^2 f(\boldsymbol{x}) +D\left(\boldsymbol{x};\mu\right)
\end{aligned}
\end{equation}

where $\boldsymbol{e}\left(\boldsymbol{x};\mu\right)=[e_i]$ and $D\left(\boldsymbol{x};\mu\right)=[d_{ij}]$ is a diagonal matrix and their entries are given by, 
\begin{equation}
\label{eq15}
e_i=\frac{1}{m q_i} \left(\frac{x_{i}-r_{i}}{q_{i}}\right)^{\mu-1} ,
d_{ii}=\frac{\mu-1}{m q^2_i} \left(\frac{x_{i}-r_{i}}{q_{i}}\right)^{\mu-2}
\end{equation}
where $r_i$ and $q_i$ are as given in (\ref{eqtx3}).

Note that the Hessian matrix is always positive definite if $f(\boldsymbol{x})$ is strictly convex or if it is convex and $x_i\neq r_i$ for all $i=1\cdots m$.

At an outer iteration $k$ and an inner iteration $j$, the Newton search direction, $\boldsymbol{p}^{k,j}$, can be obtained by solving the following linear system using Cholesky decomposition,
\begin{equation}
\label{eq16}
\nabla^2 P(\boldsymbol{x}^{k,j};\mu_k)\boldsymbol{p}^{k,j} =-\nabla P(\boldsymbol{x}^{k,j};\mu_k)
\end{equation}

One of the best methods to ensure a global convergence of Newton's method is to invoke an efficient line search \cite{Nocedal06}, in which a step length, $\alpha$, is taken along the Newton direction as,

\begin{equation}
\label{eq17}
\boldsymbol{x}^{k,j+1} =\boldsymbol{x}^{k,j}+\alpha\boldsymbol{p}^{k,j}
\end{equation}

There are many inexact line search techniques that are practically efficient of achieving this goal \cite{Nocedal06}. In general, there is a trade-off in choosing the optimal value of $\alpha$. It is desirable to choose $\alpha$ to give a substantial reduction of $P(\boldsymbol{x}^{k,j};\mu_k)$, but at the same time, this choice should be quite fast. The Armijo backtracking line search (see Algorithm \ref{algo0}) is simple but yet efficient in achieving both. It tries the unit step length $\alpha=1$ first, so it is used if it leads to a sufficient decrease in $P(\boldsymbol{x}^{k,j};\mu_k)$, otherwise it is decreased by a factor of $\rho$. Typical values used for $c$ and $\rho$ are $10^{-4}$ and $0.5$ respectively. It is often sufficient to find an approximate minimizer for (\ref{eq12}), so, the inner iterations can terminate if  one of the following criteria is met,
\begin{equation}
\label{eqstp2}
\begin{aligned}
\Vert \nabla P(\boldsymbol{x}^{k,j};\mu_k)\Vert_\infty &\le\epsilon_{gp}  \\
\vert P(\boldsymbol{x}^{k,j+1};\mu_k)-P(\boldsymbol{x}^{k,j};\mu_k)\vert &\le \epsilon_p (1+\vert P(\boldsymbol{x}^{k,j};\mu_k)\vert)  \\
\Vert \boldsymbol{x}^{k,j+1}-\boldsymbol{x}^{k,j}\Vert_\infty &\le \epsilon_x (1+\Vert \boldsymbol{x}^{k,j}\Vert_\infty)
\end{aligned}
\end{equation}

where $\epsilon_{gp}$ is a tolerance for the gradient of the barrier function, $\epsilon_p$ is a relative tolerance for the barrier function and $\epsilon_x$ is a relative tolerance for the minimizer.

\begin{algorithm}[ht]
\caption{Backtracking line search}\label{algo0}
\begin{algorithmic}[1]
\Require $\rho,c \in(0,\infty)$ 
\Ensure $\alpha$ 
\State $\alpha=1$	  
\While{$P(\boldsymbol{x}^{k,j+1};\mu_k)>P(\boldsymbol{x}^{k,j};\mu_k) +c \alpha \nabla P(\boldsymbol{x}^{k,j};\mu_k)^T\boldsymbol{p}^{k,j}$} 
	    \State $\alpha \leftarrow \rho \alpha$ 
\EndWhile
\end{algorithmic}
\end{algorithm}

The following theorems study the convergence of the MB algorithm and the number of outer iterations required to achieve a certain tolerance for the optimal function value.

\begin{theorem}[Convergence]\label{thm1}
Suppose that $\boldsymbol{z}^k$ is the minimizer of $P(\boldsymbol{z}^k;\mu_k)$ in \textnormal{(\ref{eq101})} above, and that $\mu_k \to \infty $. Then every limit point $\boldsymbol{z}$ of the sequence ${\boldsymbol{z}^k}$ is a solution of the problem 
\textnormal{(\ref{eq100})}.
\end{theorem}
\begin{proof}

Assume that $\boldsymbol{z}^*$ is the solution of (\ref{eq100}); so, $-1\le z_i^*\le1, i=1\cdots m$ and we have,

\begin{equation}
\label{eq22}
\sum_{i=1}^{m} (z_i^*)^{\mu}\le m
\end{equation}
 
Assume that $\boldsymbol{z}^k$ is the minimizer for $P(\boldsymbol{z};\mu_k)$ for each $k$, then we have, $P(\boldsymbol{z}^k;\mu_k)\le P(\boldsymbol{z}^*;\mu_k)$; that is,

\begin{equation}
\label{eq23}
f(\boldsymbol{z}^k)+ \frac{1}{m \mu_k}\sum_{i=1}^{m} (z^k_i)^{\mu_k}\le f(\boldsymbol{z}^*)+ \frac{1}{m \mu_k}\sum_{i=1}^{m} (z^*_i)^{\mu_k}
\end{equation}

Using (\ref{eq22}),
\begin{equation}
\label{eq24}
f(\boldsymbol{z}^k)+ \frac{1}{m \mu_k}\sum_{i=1}^{m} (z^k_i)^{\mu_k}\le f(\boldsymbol{z}^*)+ \frac{1}{\mu_k}
\end{equation}

Suppose that $\hat{\boldsymbol{z}}$ is a limit point for $\left\{z^k\right\}$ so that there is an infinite subsequence $\mathcal{K}$ such that,
\begin{equation}
\label{eq25}
\lim_{k \to \mathcal{K}} \boldsymbol{z}^k=\hat{\boldsymbol{z}}
\end{equation}

By taking the limit as $k \to \infty$ for $k\in \mathcal{K}$, 

\begin{equation}
\label{eq26}
f(\hat{\boldsymbol{z}})+ \lim_{k \in \mathcal{K}}\frac{1}{m \mu_k}\sum_{i=1}^{m} (\hat{z}_i)^{\mu_k}\le f(\boldsymbol{z}^*)+ \lim_{k \in \mathcal{K}}\frac{1}{\mu_k}=f(\boldsymbol{z}^*)
\end{equation}

Since the second term in the left-hand-side is always nonnegative, we have,
\begin{equation}
\label{eq27}
f(\hat{\boldsymbol{z}})\le f(\boldsymbol{z}^*)
\end{equation}
 
Thus,
\begin{equation}
\label{eq28}
f(\hat{\boldsymbol{z}})= f(\boldsymbol{z}^*)
\end{equation}

\end{proof}
\qed

\begin{theorem}[Number of outer iterations]\label{thm2}
\end{theorem}
To achieve a relative function tolerance, $\epsilon_f$, such that $\vert f(\boldsymbol{z}^n)-f(\boldsymbol{z}^*)\vert\le\epsilon_f (1+\vert f(\boldsymbol{z}^*)\vert)$ for algorithm \ref{algo1}, at least $n+1$ iterations are required where, 

\begin{equation}
\label{eq29}
n\ge\frac{log(2/(\mu_0\epsilon_f))}{log\,\tau} 
\end{equation}

\begin{proof}
From (\ref{eq24}) and for large $\mu_n$,
\begin{equation}
\label{eq30}
\vert f(\boldsymbol{z}^n)-f(\boldsymbol{z}^*)\vert \le \frac{1}{m \mu_n}\sum_{i=1}^{m}\vert (z^n_i)^{\mu_n}-(z^*_i)^{\mu_n}\vert\le \frac{2}{\mu_n}=\frac{2}{\tau^n\mu_0}
\end{equation}
Since
\begin{equation}
\label{eq31}
1+\vert f(\boldsymbol{z}^*)\vert\ge 1
\end{equation}

So, by dividing (\ref{eq30}) by (\ref{eq31}),
\begin{equation}
\label{eq32}
\frac{\vert f(\boldsymbol{z}^n)-f(\boldsymbol{z}^*)\vert}{1+\vert f(\boldsymbol{z}^*)\vert}\le \frac{2}{\tau^n\mu_0}
\end{equation}
So, to achieve a relative function tolerance, $\epsilon_f$, at least $n+1$ iterations are required where,
\begin{equation}
\label{eq33}
\frac{2}{\tau^n\mu_0}\le \epsilon_f
\end{equation}
which yields the required result (\ref{eq29}).
\end{proof}
\qed
So, for example, if $\mu_0=2^5$, $\tau=2$, then only 17 outer iterations are required to achieve a tolerance of $\epsilon_f=2^{-20}$. In practice, not all constraints are usually active and hence the number of outer iterations required is often smaller.

\section{Experimental results}\label{NE}
In this section, the proposed method is compared with the MATLAB built-in function ``fmincon" for IP and TR methods and the LBFGSB MEX wrapper for lbfgsb v3.0 FORTRAN library \cite{Becker22}. The proposed method was implemented in Matlab where SUITESPARSE \footnote{Available at: \url{https://people.engr.tamu.edu/davis/suitesparse.html}} ``cholmod2'' function \cite{Chen08} was used to find the Newton direction using Cholesky decomposition as it was faster than MATLAB built-in function ``mldivide'' with Cholesky decomposition. Since, the sparse structure of the Hessian matrix of the monomial barrier function does not change over the iterations (inner/outer), the SUITESPARSE ``analyze'' function is used once, in the first iteration only, to retrieve the best ordering that can be used to quickly update the Cholesky decomposition in the successive iterations. For each problem, the objective function, gradient and Hessian of the objective function were supplied to IP, TR, MB while the objective function and gradient were supplied to LBFGSB. All experiments were carried out using MATLAB R2021b on a Server with 2.93~GHz Intel Xeon X5679 CPU of 12 cores and 128 GB RAM under Windows 10 operating system. All methods were terminated when the projected gradient was less than $\epsilon_g=10^{-4}$ or the maximum number of outer iterations exceeded $50,000$. The parameters used for MB were $\tau=2$, $\mu_{max}=2^{40}$, $\epsilon_{gp}=10^{-4}$, $\epsilon_p=10^{-8}$, $\epsilon_x=10^{-8}$.

The experiments were performed on the convex quadratic box-constrained problems from the CUTEST~\cite{Gould15} collection. Since the trust-region reflective method requires that $u_i>l_i$ for all variables, so for a fair comparison, each variable, $x_i$, where $l_i=u_i$ was removed from the problem by setting $x_i=l_i$ and the objective function, gradient and Hessian matrix were updated accordingly. It is worth mentioning that some problems have positive semidefinite Hessian matrices, so $\eta I_m$, a small multiple of the identity matrix, was added to ensure that they become positive definite to be able to compute the Cholesky decomposition ($\eta=10^{-15}$). 

Table \ref{tab_IT1} to Table \ref{tab_IT5} show the problem name, the number of variables (N) and the number of non-zero elements of the Hessian matrix of the objective function (NNZ). They also show the the number of outer iterations, the number of active constraints and the elapsed time for each method. A constraint is considered active at a variable $x_i$ if $\vert x_i-l_i\vert< 10^{-8}(1+\vert l_i\vert)$ or $\vert u_i-x_i\vert< 10^{-8}(1+\vert u_i\vert)$. The smallest elapsed time is shown in bold font while the second smallest is underlined. Table \ref{tab_F1} to Table \ref{tab_F5} show the optimal function value and the projected gradient obtained by each method.

It can be observed that IP, TR, MB methods converged to the optimal solution achieving the same function value and the required projected gradient tolerance for all problems, however MB was the fastest one. On the other hand, LBFGSB converged to the optimal solution for all problems except DIAGPQB with $\text{N}\ge 5000$ and DIAGPQT with $\text{N}\ge 100,000$ problems (the corresponding function values and projected gradients are double underlined in Table \ref{tab_F1}). In spite of this, it was the best method (fastest method achieving the required projected gradient tolerance) for $117$ problems out of the $154$ problems of the collection while MB was the best method for $37$ problems and very close to LBFGSB in 28 problems (less than $50\%$ extra time). 

By investigating DIAGPQB, DIAGPQE and DIAGPQT problems, it was found that the unconstrained optimal solution was strictly feasible and the penalty term of MB at the starting point was very small relative to the objective function leading to an almost quadratic barrier function, so MB achieved the optimal solution in a single iteration. For CVXBQP1, the optimal solution of the monomial barrier function obtained after the first iteration was infeasible, however, its projection onto the feasible region led directly to the optimal solution of the constrained problem. 

It can also be observed that the number of outer iterations for MB is less than those for IP and TR for most problems. This is often because MB does not require the starting point to be strictly feasible as in IP and TR and allows it to be on the boundaries, which leads to a faster convergence for the problems where most/all the constraints are active at the starting point. 

It can also be noted that most of MB time is consumed in the computation of Cholesky decomposition, so speeding up this step can significantly enhance its performance. As for the future work, a low level implementation of the proposed algorithm on GPU can be explored to benefit from the GPU implementation of the SUITESPARSE ``cholmod2'' function. 

\section{Conclusions}\label{sec14}
In this article, a novel barrier function was introduced for the box-constrained convex optimization problems to convert them into unconstrained ones. The proposed barrier function has several advantages: 1) it adds a single monomial function for both the lower and upper bounds for each variable, 2) it solves a smaller linear system than that of the other barrier functions to find the Newton direction at each iteration and hence it is faster, 3) the optimal solution of the barrier function can be infeasible, however, its projection onto the feasible region can lead directly to the optimal solution of the main problem achieving a fast convergence, 4) the algorithm is very simple and can be easily implemented, 5) it can be easily extended to in include equality constraints, which has enormous applications, using Lagrangian functions. Experiments applied to the CUTEst problems showed that the MATLAB implementation of the proposed algorithm was faster than the IP and TR algorithms for all problems. Although it was slower than LBFGSB for many problems, it converged to the optimal solution for all problems unlike LBFGSB which failed to converge for some problems. A low level implementation of the MB algorithm on GPU and/or combining the proposed method with existing methods can be explored in the future to make it more competitive.   

\begin{table}[ht]
  \centering
  \caption{The number of iterations, number of active constraints and elaplsed time of the different methods.}\label{tab_IT1}
  \resizebox{0.85\textwidth}{!}{
\begin{adjustbox}{angle=90}
\begin{tabular}{lrrrrrrrrrrrrrrrr}
    \hline\noalign{\smallskip}
    \multicolumn{1}{c}{\multirow{2}[0]{*}{Problem }} & \multicolumn{1}{c}{\multirow{2}[0]{*}{N}} & {\multirow{2}[0]{*}{NNZ}} &  \multicolumn{4}{l}{Number of Outer Iterations} && \multicolumn{4}{l}{Number of Active Constraints} && \multicolumn{4}{l}{Elapsed Time} \\
    \cline{4-7}
\cline{9-12}
\cline{14-17}
     &     &    &  \multicolumn{1}{c}{LBFGSB}  & \multicolumn{1}{c}{TR} & \multicolumn{1}{c}{IP} & \multicolumn{1}{c}{MB} && \multicolumn{1}{c}{LBFGSB}& \multicolumn{1}{c}{TR} & \multicolumn{1}{c}{IP} & \multicolumn{1}{c}{MB} && \multicolumn{1}{c}{LBFGSB}& \multicolumn{1}{c}{TR} & \multicolumn{1}{c}{IP} & \multicolumn{1}{c}{MB}\\
    \noalign{\smallskip}\hline\noalign{\smallskip}
    \multicolumn{1}{l}{CVXBQP1} & 10    & 48    & 2     & 7     & 7     & 1     &       & 10    & 10    & 10    & 10    &       & \underline{0.028} & 1.206 & 1.357 & \textbf{0.017} \\
          & 50    & 328   & 2     & 7     & 7     & 1     &       & 50    & 50    & 50    & 50    &       & \underline{0.027} & 1.147 & 1.353 & \textbf{0.016} \\
          & 100   & 672   & 2     & 7     & 7     & 1     &       & 100   & 100   & 100   & 100   &       & \underline{0.025} & 1.205 & 1.353 & \textbf{0.017} \\
          & 1000  & 6968  & 2     & 7     & 7     & 1     &       & 1000  & 1000  & 1000  & 1000  &       & \underline{0.026} & 1.212 & 1.447 & \textbf{0.023} \\
          & 10000 & 69968 & 2     & 8     & 7     & 1     &       & 10000 & 10000 & 10000 & 10000 &       & \textbf{0.030} & 1.362 & 2.978 & \underline{0.063} \\
          & 100000 & 699968 & 2     & 9     & 7     & 1     &       & 100000 & 100000 & 100000 & 100000 &       & \textbf{0.069} & 2.448 & 133.306 & \underline{0.980} \\
\hline
    \end{tabular}%
 \end{adjustbox}
 }
  \label{tab:addlabel}%
\end{table}%

\begin{table}[ht]
  \centering
  \caption{The number of iterations, number of active constraints and elaplsed time of the different methods.}\label{tab_IT2}
  \resizebox{0.8\textwidth}{!}{
\begin{adjustbox}{angle=90}
    \begin{tabular}{lrrrrrrrrrrrrrrrr}
    \hline\noalign{\smallskip}
    \multicolumn{1}{c}{\multirow{2}[0]{*}{Problem }} & \multicolumn{1}{c}{\multirow{2}[0]{*}{N}} & {\multirow{2}[0]{*}{NNZ}} &  \multicolumn{4}{l}{Number of Outer Iterations} && \multicolumn{4}{l}{Number of Active Constraints} && \multicolumn{4}{l}{Elapsed Time} \\
    \cline{4-7}
\cline{9-12}
\cline{14-17}
     &     &    &  \multicolumn{1}{c}{LBFGSB}  & \multicolumn{1}{c}{TR} & \multicolumn{1}{c}{IP} & \multicolumn{1}{c}{MB} && \multicolumn{1}{c}{LBFGSB}& \multicolumn{1}{c}{TR} & \multicolumn{1}{c}{IP} & \multicolumn{1}{c}{MB} && \multicolumn{1}{c}{LBFGSB}& \multicolumn{1}{c}{TR} & \multicolumn{1}{c}{IP} & \multicolumn{1}{c}{MB}\\
    \noalign{\smallskip}\hline\noalign{\smallskip}
    \multicolumn{1}{l}{OBSTCLAE} & 16    & 12    & 2     & 5     & 8     & 1     &       & 16    & 16    & 16    & 16    &       & \underline{0.024} & 1.154 & 1.330 & \textbf{0.015} \\
          & 100   & 288   & 16    & 28    & 16    & 23    &       & 63    & 63    & 63    & 63    &       & \textbf{0.049} & 1.450 & 1.550 & \underline{0.332} \\
          & 529   & 2121  & 36    & 100   & 19    & 23    &       & 266   & 266   & 266   & 266   &       & \textbf{0.037} & 1.496 & 1.451 & \underline{0.122} \\
          & 1024  & 4380  & 47    & 115   & 22    & 23    &       & 508   & 508   & 508   & 508   &       & \textbf{0.057} & 1.731 & 1.580 & \underline{0.238} \\
          & 5625  & 26353 & 89    & 1007  & 24    & 24    &       & 2775  & 2759  & 2757  & 2771  &       & \textbf{0.149} & 32.526 & 2.935 & \underline{1.115} \\
          & 10000 & 47628 & 120   & 1510  & 24    & 24    &       & 4963  & 4952  & 4923  & 4982  &       & \textbf{0.331} & 100.558 & 4.586 & \underline{2.724} \\
\hline
    \end{tabular}%
 \end{adjustbox}
}
  \label{tab:addlabel}%
\end{table}%

\begin{table}[ht]
  \centering
  \caption{The number of iterations, number of active constraints and elaplsed time of the different methods.}\label{tab_IT3}
  \resizebox{0.85\textwidth}{!}{
\begin{adjustbox}{angle=90}
    \begin{tabular}{lrrrrrrrrrrrrrrrr}
    \hline\noalign{\smallskip}
    \multicolumn{1}{c}{\multirow{2}[0]{*}{Problem }} & \multicolumn{1}{c}{\multirow{2}[0]{*}{N}} & {\multirow{2}[0]{*}{NNZ}} &  \multicolumn{4}{l}{Number of Outer Iterations} && \multicolumn{4}{l}{Number of Active Constraints} && \multicolumn{4}{l}{Elapsed Time} \\
    \cline{4-7}
\cline{9-12}
\cline{14-17}
     &     &    &  \multicolumn{1}{c}{LBFGSB}  & \multicolumn{1}{c}{TR} & \multicolumn{1}{c}{IP} & \multicolumn{1}{c}{MB} && \multicolumn{1}{c}{LBFGSB}& \multicolumn{1}{c}{TR} & \multicolumn{1}{c}{IP} & \multicolumn{1}{c}{MB} && \multicolumn{1}{c}{LBFGSB}& \multicolumn{1}{c}{TR} & \multicolumn{1}{c}{IP} & \multicolumn{1}{c}{MB}\\
    \noalign{\smallskip}\hline\noalign{\smallskip}
    \multicolumn{1}{l}{OBSTCLBU} & 16    & 12    & 1     & 1     & 6     & 5     &       & 16    & 16    & 16    & 16    &       & \textbf{0.023} & 1.039 & 1.303 & \underline{0.032} \\
          & 100   & 288   & 10    & 11    & 11    & 5     &       & 84    & 84    & 84    & 84    &       & \textbf{0.031} & 1.187 & 1.347 & \underline{0.038} \\
          & 529   & 2121  & 21    & 10    & 15    & 7     &       & 277   & 277   & 277   & 277   &       & \textbf{0.035} & 1.193 & 1.402 & \underline{0.083} \\
          & 1024  & 4380  & 26    & 10    & 16    & 7     &       & 475   & 475   & 475   & 475   &       & \textbf{0.037} & 1.200 & 1.500 & \underline{0.151} \\
          & 5625  & 26353 & 52    & 17    & 15    & 7     &       & 1805  & 1805  & 1801  & 1805  &       & \textbf{0.097} & 1.414 & 2.539 & \underline{0.699} \\
          & 10000 & 47628 & 58    & 21    & 15    & 7     &       & 2974  & 2972  & 2966  & 2973  &       & \textbf{0.170} & 1.786 & 3.517 & \underline{0.994} \\
\hline
    \end{tabular}%
 \end{adjustbox}
}  \label{tab:addlabel}%
\end{table}%

\begin{table}[ht]
  \centering
  \caption{The number of iterations, number of active constraints and elaplsed time of the different methods.}\label{tab_IT4}
  \resizebox{0.8\textwidth}{!}{
\begin{adjustbox}{angle=90}
    \begin{tabular}{lrrrrrrrrrrrrrrrr}
    \hline\noalign{\smallskip}
    \multicolumn{1}{c}{\multirow{2}[0]{*}{Problem }} & \multicolumn{1}{c}{\multirow{2}[0]{*}{N}} & {\multirow{2}[0]{*}{NNZ}} &  \multicolumn{4}{l}{Number of Outer Iterations} && \multicolumn{4}{l}{Number of Active Constraints} && \multicolumn{4}{l}{Elapsed Time} \\
    \cline{4-7}
\cline{9-12}
\cline{14-17}
     &     &    &  \multicolumn{1}{c}{LBFGSB}  & \multicolumn{1}{c}{TR} & \multicolumn{1}{c}{IP} & \multicolumn{1}{c}{MB} && \multicolumn{1}{c}{LBFGSB}& \multicolumn{1}{c}{TR} & \multicolumn{1}{c}{IP} & \multicolumn{1}{c}{MB} && \multicolumn{1}{c}{LBFGSB}& \multicolumn{1}{c}{TR} & \multicolumn{1}{c}{IP} & \multicolumn{1}{c}{MB}\\
    \noalign{\smallskip}\hline\noalign{\smallskip}
    \multicolumn{1}{l}{TORSION4} & 16    & 12    & 2     & 4     & 5     & 1     &       & 16    & 16    & 16    & 16    &       & \textbf{0.024} & 1.112 & 1.231 & \underline{0.025} \\
          & 100   & 288   & 6     & 13    & 11    & 7     &       & 88    & 88    & 88    & 88    &       & \textbf{0.030} & 1.202 & 1.335 & \underline{0.037} \\
          & 484   & 1920  & 15    & 12    & 12    & 8     &       & 372   & 372   & 372   & 372   &       & \textbf{0.033} & 1.194 & 1.363 & \underline{0.054} \\
          & 1024  & 4380  & 21    & 14    & 12    & 8     &       & 748   & 748   & 748   & 748   &       & \textbf{0.035} & 1.216 & 1.417 & \underline{0.107} \\
          & 5476  & 25632 & 54    & 20    & 13    & 7     &       & 3692  & 3692  & 3692  & 3692  &       & \textbf{0.085} & 1.522 & 2.094 & \underline{0.404} \\
          & 10000 & 47628 & 69    & 22    & 13    & 6     &       & 6556  & 6556  & 6540  & 6564  &       & \textbf{0.164} & 1.937 & 2.808 & \underline{0.647} \\
\hline
    \end{tabular}%
 \end{adjustbox}
}
  \label{tab:addlabel}%
\end{table}%

\begin{table}[ht]
  \centering
  \caption{The number of iterations, number of active constraints and elaplsed time of the different methods.}\label{tab_IT5}
  \resizebox{0.9\textwidth}{!}{
\begin{adjustbox}{angle=90}
    \begin{tabular}{lrrrrrrrrrrrrrrrr}
    \hline\noalign{\smallskip}
    \multicolumn{1}{c}{\multirow{2}[0]{*}{Problem }} & \multicolumn{1}{c}{\multirow{2}[0]{*}{N}} & {\multirow{2}[0]{*}{NNZ}} &  \multicolumn{4}{l}{Number of Outer Iterations} && \multicolumn{4}{l}{Number of Active Constraints} && \multicolumn{4}{l}{Elapsed Time} \\
    \cline{4-7}
\cline{9-12}
\cline{14-17}
     &     &    &  \multicolumn{1}{c}{LBFGSB}  & \multicolumn{1}{c}{TR} & \multicolumn{1}{c}{IP} & \multicolumn{1}{c}{MB} && \multicolumn{1}{c}{LBFGSB}& \multicolumn{1}{c}{TR} & \multicolumn{1}{c}{IP} & \multicolumn{1}{c}{MB} && \multicolumn{1}{c}{LBFGSB}& \multicolumn{1}{c}{TR} & \multicolumn{1}{c}{IP} & \multicolumn{1}{c}{MB}\\
    \noalign{\smallskip}\hline\noalign{\smallskip}
    \multicolumn{1}{l}{TORSIONB} & 16    & 12    & 3     & 7     & 8     & 1     &       & 12    & 12    & 12    & 12    &       & \underline{0.026} & 1.163 & 1.310 & \textbf{0.014} \\
          & 100   & 288   & 9     & 12    & 12    & 8     &       & 60    & 60    & 60    & 60    &       & \textbf{0.031} & 1.189 & 1.325 & \underline{0.035} \\
          & 484   & 1920  & 26    & 11    & 12    & 8     &       & 212   & 212   & 212   & 212   &       & \textbf{0.035} & 1.196 & 1.359 & \underline{0.053} \\
          & 1024  & 4380  & 40    & 12    & 13    & 8     &       & 404   & 404   & 404   & 404   &       & \textbf{0.041} & 1.215 & 1.422 & \underline{0.087} \\
          & 5476  & 25632 & 88    & 15    & 13    & 7     &       & 1916  & 1876  & 1876  & 1884  &       & \textbf{0.132} & 1.510 & 2.089 & \underline{0.318} \\
          & 10000 & 47628 & 108   & 19    & 12    & 7     &       & 3444  & 3350  & 3308  & 3332  &       & \textbf{0.279} & 2.043 & 2.690 & \underline{0.611} \\
\hline
    \end{tabular}%
 \end{adjustbox}
}
  \label{tab:addlabel}%
\end{table}%

\begin{table}[ht]
  \centering
  \caption{The optimal function value and projected gradient of the different methods.}\label{tab_F1}
  \resizebox{0.85\textwidth}{!}{
\begin{adjustbox}{angle=90}
    \begin{tabular}{lrrrrrrrrrr}
    \hline\noalign{\smallskip}
    \multicolumn{1}{c}{\multirow{2}[0]{*}{Problem }} & \multicolumn{1}{c}{\multirow{2}[0]{*}{N}} &  \multicolumn{4}{l}{Optimal Function Value} & & \multicolumn{4}{l}{Projected Gradient } \\
\cline{3-6}
\cline{8-11}
     &       &  \multicolumn{1}{l}{LBFGSB}  & \multicolumn{1}{l}{TR} & \multicolumn{1}{l}{IP} & \multicolumn{1}{l}{MB} & & \multicolumn{1}{l}{LBFGSB}& \multicolumn{1}{l}{TR} & \multicolumn{1}{l}{IP} & \multicolumn{1}{l}{MB} \\
    \noalign{\smallskip}\hline\noalign{\smallskip}
    \multicolumn{1}{l}{CVXBQP1} & 10    & 2.4750E+00 & 2.4751E+00 & 2.4751E+00 & 2.4750E+00 &       & 0.0000E+00 & 1.2350E-05 & 8.6577E-06 & 0.0000E+00 \\
          & 50    & 5.7375E+01 & 5.7375E+01 & 5.7376E+01 & 5.7375E+01 &       & 0.0000E+00 & 3.3741E-06 & 9.5858E-06 & 0.0000E+00 \\
          & 100   & 2.2725E+02 & 2.2725E+02 & 2.2725E+02 & 2.2725E+02 &       & 0.0000E+00 & 7.0174E-06 & 9.7112E-06 & 0.0000E+00 \\
          & 1000  & 2.2523E+04 & 2.2523E+04 & 2.2523E+04 & 2.2523E+04 &       & 0.0000E+00 & 4.0059E-05 & 9.8991E-06 & 0.0000E+00 \\
          & 10000 & 2.2502E+06 & 2.2502E+06 & 2.2502E+06 & 2.2502E+06 &       & 0.0000E+00 & 1.0158E-05 & 9.9214E-06 & 0.0000E+00 \\
          & 100000 & 2.2500E+08 & 2.2503E+08 & 2.2500E+08 & 2.2500E+08 &       & 0.0000E+00 & 4.6115E-05 & 9.9237E-06 & 0.0000E+00 \\
\hline
    \end{tabular}
    \end{adjustbox}
}
  \label{tab:addlabel}%
\end{table}%

\begin{table}[ht]
  \centering
  \caption{The optimal function value and projected gradient of the different methods.}\label{tab_F2}
  \resizebox{0.8\textwidth}{!}{
\begin{adjustbox}{angle=90}
    \begin{tabular}{lrrrrrrrrrr}
    \hline\noalign{\smallskip}
    \multicolumn{1}{c}{\multirow{2}[0]{*}{Problem }} & \multicolumn{1}{c}{\multirow{2}[0]{*}{N}} &  \multicolumn{4}{l}{Optimal Function Value} & & \multicolumn{4}{l}{Projected Gradient } \\
\cline{3-6}
\cline{8-11}
     &       &  \multicolumn{1}{l}{LBFGSB}  & \multicolumn{1}{l}{TR} & \multicolumn{1}{l}{IP} & \multicolumn{1}{l}{MB} & & \multicolumn{1}{l}{LBFGSB}& \multicolumn{1}{l}{TR} & \multicolumn{1}{l}{IP} & \multicolumn{1}{l}{MB} \\
    \noalign{\smallskip}\hline\noalign{\smallskip}
    \multicolumn{1}{l}{OBSTCLAE} & 16    & 7.5366E-01 & 7.5367E-01 & 7.5367E-01 & 7.5366E-01 &       & 0.0000E+00 & 2.0424E-05 & 4.7116E-06 & 0.0000E+00 \\
          & 100   & 1.3979E+00 & 1.3979E+00 & 1.3979E+00 & 1.3979E+00 &       & 7.8193E-05 & 9.0901E-05 & 2.6964E-05 & 6.3654E-05 \\
          & 529   & 1.6780E+00 & 1.6780E+00 & 1.6780E+00 & 1.6780E+00 &       & 6.1430E-05 & 7.6159E-05 & 5.1159E-05 & 6.1439E-05 \\
          & 1024  & 1.7483E+00 & 1.7483E+00 & 1.7483E+00 & 1.7483E+00 &       & 9.1460E-05 & 9.8639E-05 & 7.7406E-05 & 6.1454E-05 \\
          & 5625  & 1.8630E+00 & 1.8630E+00 & 1.8630E+00 & 1.8630E+00 &       & 8.3798E-05 & 9.4932E-05 & 9.0148E-05 & 3.2049E-05 \\
          & 10000 & 1.8865E+00 & 1.8865E+00 & 1.8865E+00 & 1.8865E+00 &       & 8.8751E-05 & 9.9703E-05 & 9.7857E-05 & 3.1211E-05 \\
\end{tabular}
    \end{adjustbox}
}
  \label{tab:addlabel}%
\end{table}%

\begin{table}[ht]
  \centering
  \caption{The optimal function value and projected gradient of the different methods.}\label{tab_F3}
  \resizebox{0.85\textwidth}{!}{
\begin{adjustbox}{angle=90}
    \begin{tabular}{lrrrrrrrrrr}
    \hline\noalign{\smallskip}
    \multicolumn{1}{c}{\multirow{2}[0]{*}{Problem }} & \multicolumn{1}{c}{\multirow{2}[0]{*}{N}} &  \multicolumn{4}{l}{Optimal Function Value} & & \multicolumn{4}{l}{Projected Gradient } \\
\cline{3-6}
\cline{8-11}
     &       &  \multicolumn{1}{l}{LBFGSB}  & \multicolumn{1}{l}{TR} & \multicolumn{1}{l}{IP} & \multicolumn{1}{l}{MB} & & \multicolumn{1}{l}{LBFGSB}& \multicolumn{1}{l}{TR} & \multicolumn{1}{l}{IP} & \multicolumn{1}{l}{MB} \\
    \noalign{\smallskip}\hline\noalign{\smallskip}
    \multicolumn{1}{l}{OBSTCLBU} & 16    & -8.1108E-03 & -8.1108E-03 & -8.1026E-03 & -8.0908E-03 &       & 0.0000E+00 & 0.0000E+00 & 2.2488E-05 & 5.5109E-05 \\
          & 100   & 2.8750E+00 & 2.8750E+00 & 2.8751E+00 & 2.8752E+00 &       & 3.8464E-05 & 7.3274E-05 & 1.8691E-05 & 5.7530E-05 \\
          & 529   & 6.5193E+00 & 6.5193E+00 & 6.5193E+00 & 6.5193E+00 &       & 9.6207E-05 & 3.5900E-05 & 3.2829E-05 & 6.2067E-05 \\
          & 1024  & 6.8871E+00 & 6.8871E+00 & 6.8871E+00 & 6.8871E+00 &       & 9.5532E-05 & 5.1539E-05 & 4.6043E-05 & 5.9680E-05 \\
          & 5625  & 7.2309E+00 & 7.2309E+00 & 7.2309E+00 & 7.2309E+00 &       & 7.2300E-05 & 4.3083E-05 & 7.8737E-05 & 6.9110E-05 \\
          & 10000 & 7.2722E+00 & 7.2722E+00 & 7.2722E+00 & 7.2722E+00 &       & 7.9233E-05 & 5.8248E-05 & 8.6888E-05 & 7.0571E-05 \\
\hline
    \end{tabular}%
    \end{adjustbox}
}
  \label{tab:addlabel}%
\end{table}%

\begin{table}[ht]
  \centering
  \caption{The optimal function value and projected gradient of the different methods.}\label{tab_F4}%
  \resizebox{0.8\textwidth}{!}{
  \begin{adjustbox}{angle=90}
    \begin{tabular}{lrrrrrrrrrr}
    \hline\noalign{\smallskip}
    \multicolumn{1}{c}{\multirow{2}[0]{*}{Problem }} & \multicolumn{1}{c}{\multirow{2}[0]{*}{N}} &  \multicolumn{4}{l}{Optimal Function Value} & & \multicolumn{4}{l}{Projected Gradient } \\
\cline{3-6}
\cline{8-11}
     &       &  \multicolumn{1}{l}{LBFGSB}  & \multicolumn{1}{l}{TR} & \multicolumn{1}{l}{IP} & \multicolumn{1}{l}{MB} & & \multicolumn{1}{l}{LBFGSB}& \multicolumn{1}{l}{TR} & \multicolumn{1}{l}{IP} & \multicolumn{1}{l}{MB} \\
    \noalign{\smallskip}\hline\noalign{\smallskip}
    \multicolumn{1}{l}{TORSION4} & 16    & -1.2593E+00 & -1.2591E+00 & -1.2592E+00 & -1.2593E+00 &       & 0.0000E+00 & 5.1033E-05 & 1.8739E-05 & 0.0000E+00 \\
          & 100   & -1.2705E+00 & -1.2705E+00 & -1.2705E+00 & -1.2705E+00 &       & 3.3307E-16 & 6.4616E-05 & 4.1071E-05 & 6.5213E-05 \\
          & 484   & -1.2422E+00 & -1.2422E+00 & -1.2422E+00 & -1.2422E+00 &       & 9.4393E-05 & 6.5955E-05 & 7.4704E-05 & 5.6795E-05 \\
          & 1024  & -1.2317E+00 & -1.2317E+00 & -1.2317E+00 & -1.2317E+00 &       & 9.6424E-05 & 3.7433E-05 & 5.5724E-05 & 6.1682E-05 \\
          & 5476  & -1.2170E+00 & -1.2170E+00 & -1.2169E+00 & -1.2169E+00 &       & 7.5637E-05 & 6.7726E-05 & 9.6322E-05 & 7.8092E-05 \\
          & 10000 & -1.2138E+00 & -1.2138E+00 & -1.2138E+00 & -1.2138E+00 &       & 9.1174E-05 & 6.7135E-05 & 8.1052E-05 & 7.7799E-05 \\
\hline
    \end{tabular}%
    \end{adjustbox}
}
  \label{tab:addlabel}%
\end{table}%

\begin{table}[ht]
  \centering
  \caption{The optimal function value and projected gradient of the different methods.}\label{tab_F5}%
  \resizebox{0.85\textwidth}{!}{
  \begin{adjustbox}{angle=90}
    \begin{tabular}{lrrrrrrrrrr}
    \hline\noalign{\smallskip}
    \multicolumn{1}{c}{\multirow{2}[0]{*}{Problem }} & \multicolumn{1}{c}{\multirow{2}[0]{*}{N}} &  \multicolumn{4}{l}{Optimal Function Value} & & \multicolumn{4}{l}{Projected Gradient } \\
\cline{3-6}
\cline{8-11}
     &       &  \multicolumn{1}{l}{LBFGSB}  & \multicolumn{1}{l}{TR} & \multicolumn{1}{l}{IP} & \multicolumn{1}{l}{MB} & & \multicolumn{1}{l}{LBFGSB}& \multicolumn{1}{l}{TR} & \multicolumn{1}{l}{IP} & \multicolumn{1}{l}{MB} \\
    \noalign{\smallskip}\hline\noalign{\smallskip}
    \multicolumn{1}{l}{TORSIONB} & 16    & -3.0864E-01 & -3.0864E-01 & -3.0864E-01 & -3.0864E-01 &       & 1.1102E-16 & 9.5096E-07 & 3.3215E-05 & 3.3307E-16 \\
          & 100   & -4.0570E-01 & -4.0570E-01 & -4.0570E-01 & -4.0569E-01 &       & 9.1601E-05 & 4.9096E-05 & 2.2632E-05 & 8.7275E-05 \\
          & 484   & -4.1611E-01 & -4.1611E-01 & -4.1611E-01 & -4.1610E-01 &       & 8.4080E-05 & 7.5153E-05 & 8.4724E-05 & 6.2527E-05 \\
          & 1024  & -4.1740E-01 & -4.1740E-01 & -4.1739E-01 & -4.1738E-01 &       & 8.1777E-05 & 8.8648E-05 & 6.8616E-05 & 6.3897E-05 \\
          & 5476  & -4.1829E-01 & -4.1829E-01 & -4.1829E-01 & -4.1827E-01 &       & 9.8695E-05 & 8.6701E-05 & 8.5754E-05 & 9.3502E-05 \\
          & 10000 & -4.1834E-01 & -4.1839E-01 & -4.1837E-01 & -4.1837E-01 &       & 6.0746E-05 & 9.5136E-05 & 9.2585E-05 & 7.4569E-05 \\
\hline
    \end{tabular}%
    \end{adjustbox}
}
  \label{tab:addlabel}%
\end{table}%

\bibliographystyle{unsrt}
\bibliography{mb_arxiv}  

\end{document}